\documentclass[reqno, letterpaper]{amsart}
\usepackage{amssymb}
\usepackage{amsthm}
\usepackage{hyperref}
\usepackage{cleveref}
\hypersetup{
	colorlinks = true,
	linkcolor = blue,
	citecolor = magenta
}

\usepackage[OT2,T1]{fontenc}
\DeclareSymbolFont{cyrletters}{OT2}{wncyr}{m}{n}
\DeclareMathSymbol{\Sha}{\mathalpha}{cyrletters}{"58}

\usepackage{tikz-cd}
\tikzcdset{
	diagrams={sep=small}			% small default diagram size
}

\usepackage{multirow}				% merging cells vertically
\usepackage{mathtools}				% more math symobls

\usepackage{enumitem}
\setlist[1]{labelindent=\parindent}
\setlist[enumerate, 1]{label = \textnormal{(\arabic*)}, ref = \arabic*}
\setlist[enumerate, 2]{label = \textnormal{(\alph*)}, ref = \alph*}
\setlist[enumerate, 3]{label = \textnormal{(\roman*)}, ref = \roman*}

\makeatletter
	\def\l@subsection{\@tocline{2}{0pt}{2.5pc}{5pc}{}}
\makeatother

\theoremstyle{plain}
\newtheorem{theorem}[equation]{Theorem}
\newtheorem{proposition}[equation]{Proposition}
\newtheorem{lemma}[equation]{Lemma}
\newtheorem{corollary}[equation]{Corollary}

\theoremstyle{definition}
\newtheorem{definition}[equation]{Definition}

\theoremstyle{remark}
\newtheorem{remark}[equation]{Remark}

\numberwithin{equation}{section}

\let\emptyset\varnothing

\let\AA\undefined

\newcommand{\AA}{\mathbb{A}}
\newcommand{\CC}{\mathbb{C}}

\newcommand{\PP}{\mathbb{P}}
\newcommand{\QQ}{\mathbb{Q}}

\DeclareMathOperator{\Aut}{Aut}

\DeclareMathOperator{\id}{id}
\DeclareMathOperator{\im}{im}
\DeclareMathOperator{\coker}{coker}
\DeclareMathOperator{\rk}{rk}

\DeclareMathOperator{\Sym}{Sym}

\DeclareMathOperator{\Spec}{Spec}

\DeclareMathOperator{\Pic}{Pic}

\DeclareMathOperator{\Alb}{Alb}

\DeclareMathOperator{\CH}{CH}

\DeclareMathOperator{\ch}{ch}

\DeclareMathOperator{\alb}{alb}

\usepackage[centering]{geometry}	% page layout and margins

\usepackage[style=alphabetic, doi=false, isbn=false, url=false]{biblatex}
\DeclareMathOperator{\forget}{forget}

\begin{document}
\title[OpenAI's base of Lagrangian fibrations]{Note on OpenAI's proof of the base of Lagrangian fibrations}

\author[Y.-J. Kim]{Yoon-Joo Kim}
\address{Department of Mathematics, Columbia University, New York, NY 10027, USA}
\address{HCMC, Korea Institute for Advanced Study, Seoul, South Korea}
\email{yk3029@columbia.edu}

%\date{\today}

\begin{abstract}
	OpenAI \cite{OAI:Pn} proved that every Lagrangian fibration of a compact hyper-K\"ahler manifold has a base $\PP^n$. This expository note summarizes the last section of their paper.
\end{abstract}

\maketitle
%\tableofcontents

\subsection*{Disclosures}
	I declare that all work here is my human effort on understanding the paper \cite{OAI:Pn}. No AI was used in this process, neither to understand the paper nor to write up this note (except that \cite{OAI:Pn} was entirely written by AI). All the results here are either explicitly or tacitly contained in OpenAI's paper, and I claim no originality. Consider this as an expository seminar note.
	
	Shortly after OpenAI's paper has appeared, Zhiyu Liu announced a separate paper \cite{liu26}. Their Theorem~1.1 extracted an optimal statement from OpenAI's argument. I will follow Liu's formulation, but will exclusively follow OpenAI's paper for the body of the note.
%	Their Theorem~1.1 extracted an optimal statement from OpenAI's argument, and the main statement of this note almost inevitably converged to theirs. I will thus follow the nice formulation of \cite[Theorem~1.1]{liu26}. However, let me stress that no proofs in this note are influenced by Liu's result; all results are contained in the original paper \cite{OAI:Pn}.

\section{Introduction}
	The goal of this note is to understand the proof of \cite[Proposition~11.1]{OAI:Pn}. We will reference this paper many times throughout the note, so we will simply write it as AI::Proposition~11.1.
	
	\begin{theorem} [AI::Proposition~11.1] \label{thm:base}
		Let $f : X \to B$ be a Lagrangian fibration of a projective hyper-K\"ahler manifold $X$ of dimension $2n$. Suppose that $B$ has quotient singularities and there exists a finite surjective morphism $p : A \to B$ from an abelian variety $A$. Then $B$ is smooth.
	\end{theorem}
	
	This theorem is merely the last step of OpenAI's main result, where the singularity assumption on $B$ and the existence of $p$ was completely dropped. More than half of the paper is devoted to proving $B$ has quotient singularities, and AI::Proposition~9.1 in the paper shows such a morphism $p : A \to B$ always exists.
%	This part depends on OpenAI's another paper proving the hyper-K\"ahler abundance conjecture, which was proved separately in \cite{engel-mauri26}.
	We will not discuss these parts at all, and will only focus on proving their AI::Proposition~11.1 (\Cref{thm:base} in this note).
	
	\bigskip
	
	The first step of the proof is to take the canonical smooth DM stack $\mathfrak B \to B$ and lift $f$ and $p$ to $\tilde f : X \to \mathfrak B$ and $\tilde p : A \to \mathfrak B$ (AI::Lemma~11.2) using \cite[Theorem~4.6]{fan-mann-nir10}. Note that $\tilde f$ and $\tilde p$ are flat because $\mathfrak B$ is smooth. To ease notation, we will write $B$ for the canonical smooth DM stack from now on, and write
	\[ f : X \to B ,\qquad p : A \to B \]
	for the new stacky-versions of these morphisms. The original base variety of the Lagrangian fibration in \Cref{thm:base} will be denoted by $B^{sch}$, the coarse moduli space of $B$.
	
	\Cref{thm:base} rarely uses any input from hyper-K\"ahler manifolds or Lagrangian fibrations. In fact, a more technically precise version of OpenAI's result is the following, which was already extracted in \cite[Theorem~1.1]{liu26}.
	
	\begin{theorem} \label{thm:main}
		Let $B$ be a smooth proper DM stack with trivial stabilizers in codimension $1$, and admitting a finite flat morphism $p : A \to B$ from an abelian variety $A$. Assume
		\begin{enumerate}[label=\textnormal{(\roman*)}]
			\item The coarse moduli space $B^{sch}$ is a Fano variety.
			\item Every connected component of the inertia stack $I_B$, except the identity component, has the Euler characteristic $0$.
		\end{enumerate}
		Then $B$ is a smooth projective variety.
	\end{theorem}
	
	Here the Euler characteristic of a DM stack means a stacky-Euler characteristic that can be a rational number in general. The stacky Gauss--Bonnet theorem \cite[Proposition~1.6]{beh09} shows that (ii) is equivalent to the vanishing of $c_{top}(T_{I_B})$ on each non-identity component of $I_B$. Hence \Cref{thm:main} is identical to \cite[Theorem~1.1]{liu26}. OpenAI's proof of this theorem is extremely systematic and geometric. However, the paper is difficult to grasp at first sight due to its terse presentation. Throughout the note, we will expand OpenAI's original proof and make its beautiful geometric ideas more explicit.
	
	\medskip
	
	Here is a quick summary of the structure of this note. \S \ref{sec:basics} reviews and introduces basic definitions we use throughout the note: three keywords are groupoid schemes, their inertial components, and the difference morphisms. In the following two sections, we prove \Cref{thm:delta(R_i)=A} assuming condition (i) in \Cref{thm:main}, and prove \Cref{thm:delta(R_i) is not A} from condition (ii). If $B$ is not a variety then these two conclusions are incompatible, thereby completing the proof of the theorem. The last \S \ref{sec:proof of main theorems} proves the desired \Cref{thm:base}.

	\subsection*{Acknowledgments}
		I thank Andres Fernandez Herrero, James Hotchkiss, and Niklas M\"uller for helpful comments on the draft of this note.

	\subsection*{Notation and conventions}
		We work over $\CC$. A \emph{variety} in this note is a reduced and separated scheme of finite type over $\CC$; it may be reducible or disconnected unless stated otherwise.
		
		We mostly follow OpenAI's original notation, so that the reader can easily compare our exposition to the original reference. Some of the original notation looked a bit awkward to us, and in those cases, we tried to choose the most reasonable alternative.

\section{Inertial components and the difference map} \label{sec:basics}
	This section focuses on the abstract study of finite flat morphisms $p : A \to B$ from an abelian variety $A$ to a smooth proper DM stack $B$ (with no stabilizers in codimension $1$). No additional assumptions will be imposed in this section.
	
	We will define the three main objects of this note. First, we will recall the definition of the \emph{groupoid scheme} $R \rightrightarrows A$ associated to the morphism $p$. Second, we will introduce the notion of its \emph{inertial component} $R_i$ (\Cref{def:components of R}). Third, we will define the \emph{difference morphism} $\delta : R \to A$ (\Cref{def:difference}). The proof of \Cref{thm:main} in one sentence will be:
	\begin{center}
		the study of inertial components $R_i \subset R$ under the difference morphism $\delta : R \to A$.
	\end{center}

	\subsection{Finite flat groupoid schemes}
		We will heavily use the theory of algebraic stacks, or more precisely the groupoid-quotient aspects of them. Our main reference is the chapters ``Groupoid schemes'' and ``Groupoids in algebraic spaces'' in \cite[Tag~022L, 0437]{stacks-project}, which in turn references the original paper \cite{keel-mori97}.
		
		Observe that the finite flat morphism $p : A \to B$ can be thought of as a finite flat presentation of the DM stack $B$. This allows one to view $B$ as a flat quotient stack $B = [A/R]$ of the abelian variety $A$ by a \emph{groupoid scheme} $R$ over it. The latter is defined as the fiber product
		\[\begin{tikzcd}
			R \coloneq A \times_B A \arrow[r, "t"] \arrow[d, "s"] & A \arrow[d, "p"] \\
			A \arrow[r, "p"] & B
		\end{tikzcd}.\]
		The two projection maps $s, t : R \to A$ are finite flat morphisms of schemes. Take the irreducible decomposition
		\[ R = \bigcup_i R_i \qquad\text{(with reduced induced structures)} .\]
		The irreducible components $R_i$ will play very important roles throughout.
		
		\begin{lemma} \label{lem:R}
			$R$ is a projective Cohen--Macaulay variety with pure dimension $n$. Moreover, for every irreducible component $R_i \subset R$, the restricted morphisms $s_{|R_i}, t_{|R_i} : R_i \to A$ are finite surjective.
		\end{lemma}
		\begin{proof}
			Since $R$ is finite flat over $A$, it is a projective Cohen--Macaulay scheme and all its irreducible components $R_i$ are finite surjective over $A$. Since $p : A \to B$ is generically \'etale, $R$ satisfies the condition $(\mathrm R_0)$ and hence it is reduced.
		\end{proof}
		
		The variety $R$ parametrizes triples $(x,y,g)$ of the form
		\[ R = \{ (x,y,g) : \ x,y \in A ,\ \ p(x) = p(y) \in B ,\ \ g \in \Aut(p(x)) \} .\]
%		$(x,y,g)$ where $x, y \in A$ are points with $p(x) = p(y) \in B$, and $g$ is an automorphism of this stacky-point.
		The forgetful map defines a finite morphism
		\[ \forget \coloneq (s,t) : R \to A \times A ,\qquad (x,y,g) \mapsto (x,y) .\]
		We denote the image of $R$ in $A \times A$ by $\bar R = \forget(R)$, as well as $\bar R_i = \forget(R_i)$.
		
		\begin{lemma} \label{lem:bar R_i}
			The forgetful morphism $R \to A \times A$ is finite and generically injective. In particular, $R_i \to \bar R_i$ is a finite birational morphism for each $i$, so they share the same normalization $Z_i$.\footnote{The original paper often uses the normalization $Z_i$ instead of $R_i$. We will try to stick to the original variety $R_i$ as much as possible.}
		\end{lemma}
		\begin{proof}
			The DM stack $B$ has a trivial generic stabilizer, so $R \to A \times A$ is generically injective on each irreducible component $R_i$.
		\end{proof}

		\subsubsection{Compositions in groupoid}
			The groupoid $R$ comes equipped with a composition morphism in its datum. More precisely, consider the fiber product
			\[\begin{tikzcd}
				R \times_A R \arrow[r, "p_2"] \arrow[d, "p_1"] & R \arrow[d, "s"] \\
				R \arrow[r, "t"] & A
			\end{tikzcd},\qquad\quad
			\begin{aligned}
				R \times_A R &\coloneq \{ \Big( (x,y,g), (y',z,h) \Big) \in R \times R : y = y' \} \\
				&\cong A \times_B A \times_B A
			\end{aligned}.\]
			Projection map to the first-third factors is the \emph{composition} morphism
			\[ c : R \times_A R \to R ,\qquad (x,y,z, g, h) \mapsto (x,z, h \circ g) .\]
			
			\begin{lemma}
				$R \times_A R$ is a projective Cohen--Macaulay variety, whose irreducible components are $n$-dimensional and finite surjective over $A$.
			\end{lemma}
			\begin{proof}
				Same as \Cref{lem:R}.
			\end{proof}
			
			Here is a slightly stronger version:
			
			\begin{lemma} \label{lem:irreducible components of R^2}
				Let $R_i, R_j \subset R$ be irreducible components. Then the closed subscheme $R_i \times_A R_j \subset R \times_A R$ satisfies the following.\footnote{This is tacitly proved in p.94 of the paper.}
				\begin{enumerate}
					\item Every irreducible component $T \subset R_i \times_A R_j$ is an irreducible component of $R \times_A R$.
					\item $R_i = p_1(T)$, $R_j = p_2(T)$, and $R_k \coloneq c(T)$ are all irreducible components of $R$.
					\item $T$ is finite surjective over $R_i$, $R_j$, and $R_k$.
				\end{enumerate}
			\end{lemma}
			\begin{proof}
				(1) Notice that $R_i \times_A R_j$ is the preimage of $\Delta_A \subset A \times A$ under the morphism $(t,s) : R_i \times R_j \to A \times A$. Since $R_i \times R_j$ has pure dimension $2n$ and $\Delta_A \subset A \times A$ is a regular closed immersion of codimension $n$, we conclude that $R_i \times_A R_j$ has pure dimension $n$. Therefore, $R_i \times_A R_j$ is a union of irreducible components of $R \times_A R$.
				
				(2) Note that $c$ is finite, because both $R \times_A R$ and $R$ are finite over $A$. Hence $c(T)$ has dimension $n$, meaning it is an irreducible component of $R$.
				
				(3) All morphisms $s$, $t$, and $c$ are finite, so $T$ is finite over its images.
			\end{proof}
			
%			\begin{lemma}
%				Let $T \subset R \times_A R$ be an irreducible component. Then
%				\[ R_i \coloneq p_1(T) ,\quad R_j \coloneq p_2(T) ,\quad R_k \coloneq c(T) \]
%				are all irreducible components of $R$. Moreover, $T$ is finite surjective over them.
%			\end{lemma}
%			\begin{proof}
%				By the previous lemma, $T$ is an $n$-dimensional irreducible variety. Since $s$ is finite, its image $R_i = s(T) \subset R$ is again $n$-dimensional irreducible whence an irreducible component. Same holds for $R_j$ and $R_k$: note that $c$ is finite because both $R \times_A R$ and $R$ are finite over $A$.
%			\end{proof}

	\subsection{Inertial components}
		The DM stack $B$ has an inertia group stack $I_B \to B$. Its pullback
		\begin{equation} \label{diag:inertia group scheme}
		\begin{tikzcd}
			I_A \arrow[r] \arrow[d] & I_B \arrow[d] \\
			A \arrow[r, "p"] & B
		\end{tikzcd}
		\end{equation}
		is a finite group scheme $I_A \to A$ encoding every inertial information of $B$.
		
		\begin{definition}
			The finite group scheme $I_A \to A$ in \eqref{diag:inertia group scheme} is called the \emph{stabilizer group scheme} (or \emph{inertia group scheme}).\footnote{The paper calls $I_A$ the difference-zero set and did not assign it a symbol. See \Cref{lem:difference zero}.}
		\end{definition}
		
		Note that $I_A$ parametrizes the pairs $(x, g)$ where $x \in A$ and $g$ is an automorphism of the point $p(x) \in B$. Note that $I_A$ may be a non-reduced scheme. We will often consider $I_A$ as a closed subscheme of $R$ via $I_A = \forget^{-1}(\Delta_A)$, or equivalently
		\[ I_A \hookrightarrow R ,\qquad (x,g) \mapsto (x,x,g) .\]
		Under this identification, the structural morphism $I_A \to A$ is the restriction of the morphism $s : R \to A$ (or $t$)
		\[ s : I_A \to A ,\qquad (x,x,g) \mapsto x .\]
		
		\begin{definition} \label{def:components of R}
			\begin{enumerate}
				\item The \emph{identity component} of $R$ is the irreducible component
				\[ R_0 = \{ (x,x,\id) \in R \} \ \ \subset \ R \]
				that corresponds to the identity section of the group scheme $I_A \to A$.
				
				\item A \emph{non-identity component} of $R$ is an irreducible component $R_i \subset R$ that is not $R_0$.
				
				\item An \emph{inertial component} of $R$ is an irreducible component $R_i \subset R$ containing a non-identity element of $I_A$, i.e., $R_i$ contains a point of the form $(x,x,g)$ with $g \neq \id$.\footnote{The paper calls a point $(x,x,g)$ with $g \neq \id$ a ``non-identity arrow of difference zero''.}
			\end{enumerate}
		\end{definition}
		
		\begin{lemma} \label{lem:I_A}
			\begin{enumerate}
				\item The morphism $I_A \to A$ is unramified.
				\item Every connected component $S$ of $I_A$ satisfies either $S = R_0$ or $S \cap R_0 = \emptyset$.
			\end{enumerate}
		\end{lemma}
		\begin{proof}
			(1) Every finite group scheme in characteristic $0$ is unramified by Cartier's theorem.
			
			(2) Since $I_A \to A$ is unramified, \'etale locally on $I_A$, it is a closed immersion. Hence the identity component can $R_0$ never intersect with a non-identity component $S$.
		\end{proof}
		
		\begin{remark}
			A non-identity component $R_i \subset R$ may very well satisfy $R_i \cap R_0 \neq \emptyset$, see \Cref{lem:codimension 1 description of R}. This will be crucial in \Cref{lem:prop:inertial component 1}.
		\end{remark}
		
		\begin{proposition} \label{prop:inertial component}
			The following are equivalent.
			\begin{enumerate}
				\item $B \to B^{sch}$ is an isomorphism, i.e., $B$ is an algebraic space.
				\item $\forget : R \to A \times A$ is a closed immersion.
				\item The stabilizer group scheme $I_A \to A$ is trivial.
				\item $R$ has no inertial components.
			\end{enumerate}
		\end{proposition}
		\begin{proof}
			The equivalence of (1--3) and the implication (3) $\Longrightarrow$ (4) are clear. For (4) $\Longrightarrow$ (3), assume that $I_A \to A$ is nontrivial, so that there exists an irreducible component $S (\neq R_0) \subset I_A$. Since $S \subset I_A \subset R$, there exists an irreducible component $R_i \subset R$ containing $S$. Notice that $S \cap R_0 = \emptyset$ by \Cref{lem:I_A}, so $R_i$ is inertial.
		\end{proof}
		
		The following two lemmas describe the local behaviors of $R$ around the points $(x,x,g) \in I_A \subset R$. These are completely formal.
		
		\begin{lemma} \label{lem:codimension 1 description of R}
			\begin{enumerate}
				\item $R_0 \to \bar R_0 = \Delta_A$ is an isomorphism.
				\item Let $\eta \in A$ be a codimension $1$ point. Then the point $(\eta, \eta, \id) \in R$ admits an \'etale local neighborhood $W$ of the form
				\[ W \subset \{ (x,y) \in \AA^2 : x^e = y^e \} \times \AA^{n-1} ,\]
				where the identity branch $R_0 \subset R$ corresponds to $\{ (x,x) \} \times \AA^{n-1} \subset W$.
%					\[ W \subset (x^e = y^e) \times \AA^{n-1}_u ,\quad U \subset \AA^1_x \times \AA^{n-1}_u ,\qquad s(x,y,u) = (x,u) ,\quad t(x,y,u) = (y,u) .\]
			\end{enumerate}
		\end{lemma}
		\begin{proof}
			(1) Clear. (2) Locally, $p : A \to B$ is a finite morphism between smooth varieties. Therefore, around a codimension $1$ point $\eta \in A$ it has a form
			\[ p : \AA^1_x \times \AA^{n-1}_u \to \AA^1_z \times \AA^{n-1}_u ,\qquad (x,u) \mapsto (x^e, u) .\]
			The fiber product $R = A \times_B A$ around $(\eta, \eta, \id)$ therefore has the desired description.
		\end{proof}
		
		\begin{lemma} \label{lem:neighborhood of inertial component}
			Every point of the form $(x_0,x_0,g_0) \in R$ with $g_0 \neq \id$ admits an \'etale local neighborhood $W$ satisfying the following.
			\begin{enumerate}[label=\textnormal{(\alph*)}]
				\item There exists an \'etale neighborhood $U$ of $x_0 \in A$, a neighborhood $[V/G]$ of $p(x_0) \in B$, and a morphism $f : U \to V$ such that $p$ has a local presentation
				\[\begin{tikzcd}[column sep=normal]
					\bar f : U \arrow[r, "f"] & V \arrow[r] & {[V/G]}
				\end{tikzcd}.\]
				\item There exists an \'etale neighborhood $W$ of $(x_0, x_0, g_0) \in R_i$ such that every point $(x,y,g_0)$ lying in $W$ satisfies
				\[ f(y) = g_0 \cdot f(x) \qquad\text{in } \ V .\]
			\end{enumerate}
		\end{lemma}
		\begin{proof}
			The DM stack $B$ admits a neighborhood $[V/G]$ around $p(x_0)$, where $G$ is the automorphism group of $p(x_0)$ and $V$ is a smooth affine variety with a $G$-action fixing $0 \in V$ corresponding to $p(x_0)$ \cite[Lemma~2.2.3]{abr-vis02} \cite[Theorem~2.12]{ols06}. Take any neighborhood $U$ of $x \in A$ so that we have a restriction $\bar f \coloneq p_{|U} : U \to [V/G]$.
			\[\begin{tikzcd}
				& A \arrow[r, "p"] & B \\
				& U \arrow[r, "\bar f"] \arrow[u, "\acute et"] & {[V/G]} \arrow[u, "\acute et"'] \\
				U \arrow[r, hook] \arrow[ru, "\id"] & U \times G \arrow[r] \arrow[u] & V \arrow[u, "G"']
			\end{tikzcd}\]
			Shrinking $U$ further, we may assume that $V \to [V/G]$ pullbacks to $U \times G \to U$. Taking a component $U \hookrightarrow U \times G$ by $u \mapsto (u, \id)$, the third row is the desired morphism $f : U \to V$.\footnote{The paper used the notation $\psi$ (with analytic neighborhoods) in AI::Lemma~11.5.}
			
			\medskip
			
			The point $(x_0,x_0,g_0) \in R$ admits a neighborhood $R_U$ given by
			\[\begin{tikzcd}
				R_U \arrow[r, "t"] \arrow[d, "s"] & U \arrow[d, "\bar f"] \\
				U \arrow[r, "\bar f"] & {[V/G]}
			\end{tikzcd},\qquad \begin{aligned}
				R_U &= \{ (x,y,g) \in U \times U \times G : f(y) = g \cdot f(x) \} \\
				&= \bigsqcup_{g \in G} f^{-1} (V^g)
			\end{aligned}.\]
			The component $f^{-1}(V^{g_0}) \subset R_U$ is the desired neighborhood $W$ consisting of points $(x,y,g_0)$.
		\end{proof}

	\subsection{The difference morphism}
		The discussions so far were purely formal about the finite flat presentation $p : A \to B$. Let us now use the fact that $A$ is an abelian variety.
		
		\begin{definition} \label{def:difference}
			Write $A_\cdot \coloneq (A, 0)$. The \emph{difference} morphism $\delta : R \to A_\cdot$ is the composition
			\[\begin{tikzcd}[sep = normal]
				\delta : R \arrow[r, "{\forget}"] & A \times A \arrow[r, "-"] & A_\cdot ,
			\end{tikzcd} \qquad (x,y,g) \mapsto y - x .\]
			It is neither finite nor flat. We denote the restriction of the morphism $- : A \times A \to A_\cdot$ to $\bar R$ by
			\[ \bar \delta : \bar R \to A_\cdot .\]
		\end{definition}
		
		\begin{remark}
			The variety $A$ does not come equipped with a preferred choice of zero, so it has no canonical addition operation (i.e., $A$ is an abelian torsor rather than an abelian variety). The notation $A_\cdot = (A,0)$ emphasizes that $A$ is equipped with the choice of a canonical zero, so that we can consider its addition morphism $+ : A_\cdot \times A_\cdot \to A_\cdot$. Notice that the subtraction and translation morphisms
			\[ - : A \times A \to A_\cdot ,\qquad + : A_\cdot \times A \to A \]
			are always well-defined even without fixing a zero of $A$.
			
			This somewhat pedantic and annoying notation will be justified later. It will be important for us to distinguish the role of the targets of $s, t : R \to A$ and $\delta : R \to A_\cdot$. Such notation $A_\cdot$ will visually tell us the space we are working on. The discussion will be quite confusing without it.
		\end{remark}
		
		\begin{lemma} \label{lem:difference zero}
			$I_A = \delta^{-1}(0)$.
		\end{lemma}
		\begin{proof}
			Recall the description $I_A = \{ (x,x,g) \in R \}$ in $R$.
		\end{proof}
		
		\begin{proposition} \label{prop:composition}
			There exists a commutative diagram\footnote{This is tacitly mentioned in p.94 of the paper.}
			\[\begin{tikzcd}[column sep=normal]
				R \times_A R \arrow[r, "c"] \arrow[d, "\delta \times \delta"] & R \arrow[d, "\delta"] \\
				A_\cdot \times A_\cdot \arrow[r, "+"] & A_\cdot
			\end{tikzcd}.\]
			In other words, $\delta \circ c = \delta \circ p_1 + \delta \circ p_2$.
		\end{proposition}
		\begin{proof}
			This follows from the direct computation
			\begin{align*}
				(\delta \circ c) \Big( (x,y,g), (y,z,h) \Big) = \delta (x,z,h \circ g) &= z - x \\
				&= (y-x) + (z-y) = \delta(x,y,g) + \delta(y,z,h) .\qedhere
			\end{align*}
		\end{proof}

\section{Fano property of the coarse moduli space}
	In this section, we study the difference map $\delta : R \to A_\cdot$ under the condition in \Cref{thm:main}(i). Recall from \Cref{prop:inertial component} that if $B$ is not a scheme, then there exists an inertial component $R_i \subset R$. The main result of this section is a stronger version of this.
	
	\begin{theorem} \label{thm:delta(R_i)=A}
		Let $p : A \to B$ be as usual and assume \Cref{thm:main}(i). If $B$ is not a scheme, then there exists an inertial component $R_i \subset R$ with $\delta(R_i) = A_\cdot$.
	\end{theorem}
	
	Given an inertial component $R_i$, OpenAI's idea is to find a new inertial component $R_k$ with $\delta(R_i) \subsetneq \delta(R_k)$ and repeat this until $\delta(R_i) = A_\cdot$. The new inertial component $R_k$ will be a groupoid-composition of $R_i$ with an auxiliary irreducible component $R_j$. The auxiliary component $R_j$ will be found nearby the identity component $R_0$ using local geometry.

	\subsection{Image of the difference map}
		\Cref{thm:main}(i) assumes $B^{sch}$ is a Fano variety, and in particular projective. In this subsection, we will see the projectiveness of $B^{sch}$ already deduces an interesting property of $\delta : R \to A_\cdot$.
		
		\medskip
		
		Fix any ample line bundle on $B^{sch}$ and consider its pullback to $A$ via the finite morphism $A \to B \to B^{sch}$. This yields an ample line bundle $P$ on $A$ with the condition
		\begin{equation} \label{eq:sP = tP}
			s^* P = t^* P \quad\text{as line bundles on } R .
		\end{equation}
		Recall that the \emph{polarization} (with respect to $P$) is a homomorphism of abelian varieties
		\[ \lambda : A_\cdot \to \check A_\cdot ,\qquad a \mapsto [t_a^* P \otimes P^{-1}] ,\]
		where $t_a : A \to A$ is a translation automorphism by $a \in A_\cdot$, and $\check A_\cdot = \Pic^\circ_A$ is the neutral Picard scheme of $A$. Given an abelian subvariety $D \subset A_\cdot$, its \emph{complementary abelian subvariety} $D^\perp \subset A_\cdot$ is
		\[ D^\perp \coloneq \ker \big( A_\cdot \xlongrightarrow{\lambda} \check A_\cdot \twoheadrightarrow \check D \big)^\circ ,\]
		where the superscript $-^\circ$ means taking the identity component of an algebraic group. Notice the equality $(D^\perp)^\perp = D$.
		
		\medskip
		
		The following is the main result of this section.
		
		\begin{proposition} [AI::Lemma~11.6] \label{prop:image of R_i}
			Let $p : A \to B$ be as usual and assume $B^{sch}$ is projective. Then every irreducible component $R_i \subset R$ satisfies the following.
			\begin{enumerate}
				\item The closed subvariety $\delta(R_i) \subset A_\cdot$ is of the form
				\[ \delta(R_i) = D_i + \tau_i ,\]
				where $D_i \subset A_\cdot$ is an abelian subvariety and $\tau_i \in A_\cdot$ is a torsion point.
				
				\item Set $K_i = D_i^\perp$ and define a free $K_i$-action on $A \times A$ by
				\[ k \cdot (x,y) = (k + x, k + y) .\]
				Then $\bar R_i  = \forget(R_i) \subset A \times A$ is closed under the $K_i$-action.
			\end{enumerate}
		\end{proposition}
		
		\begin{corollary} \label{cor:image of R_i}
			In \Cref{prop:image of R_i}, assume further either of the following:
			\begin{enumerate}[label=\textnormal{(\roman*)}]
				\item $R_i$ is an inertial component.
				\item $R_i \cap R_0 \neq \emptyset$.
			\end{enumerate}
			Then $\delta(R_i) \subset A_\cdot$ is an abelian subvariety.
		\end{corollary}
		\begin{proof}
			The condition (i) (resp. (ii)) implies $R_i \cap I_A \neq \emptyset$, or equivalently $0 \in \delta(R_i)$. Hence the claim follows from \Cref{prop:image of R_i}(1).
		\end{proof}
		
		The proof of the proposition depends on the following key lemma.\footnote{This is extracted from AI::Lemma~11.6.}
		
		\begin{lemma} \label{lem:abelian variety}
			Let $s : F \to A$ be a morphism from a proper connected variety $F$ to an abelian variety $A$. Then
			\begin{enumerate}
				\item The pullback homomorphism $s^* : \check A_\cdot \to \Pic_F$ satisfies
				\[ \dim \ker s^* \le \dim A - \dim s(F) .\]
				Equality holds iff $s(F) \subset A$ is an abelian subtorsor (translation of an abelian subvariety).
				
				\item The same inequality holds for $\ker \big( s^* \circ \lambda : A_\cdot \to \Pic_F \big)$, and the equality holds iff
				\[ s(F) = K + \sigma \qquad\text{for}\quad K = \Big( \ker (s^* \circ \lambda)^\circ\Big)^\perp \subset A_\cdot , \ \ \sigma \in A .\]
			\end{enumerate}
		\end{lemma}
		\begin{proof}
			(1) A proper connected variety over an algebraically closed field of characteristic $0$ admits the notion of the Albanese variety \cite{con:albanese, lau-sch24}. This factorizes $s$ into
			\[\begin{tikzcd}
				F \arrow[rrr, "s"] \arrow[rd, "\alb"'] & & & A \\
				& \Alb_F \arrow[r, "s'", twoheadrightarrow] & H \arrow[ru, hook, "i"']
			\end{tikzcd},\]
			where $H \subset A$ is the image of $\Alb_F \to A$, an abelian subtorsor. Since $\alb^*$ is injective on neutral Picard schemes \cite[\S 8]{lau-sch24}, we have
			\begin{equation} \label{eq:lem:abelian variety}
				(\ker s^*)^\circ = \ker (i^* : \check A_\cdot \twoheadrightarrow \check H_\cdot) .
			\end{equation}
			This shows the dimension estimate $\dim \ker s^* = \dim A - \dim H \le \dim A - \dim s(F)$. When the equality holds, this mans $s(F) = H$ is an abelian subtorsor.
			
			\medskip
			
			(2) The dimension inequality is immediate because $\lambda : A_\cdot \to \check A_\cdot$ is an isogeny. When the equality holds, we have $s(F) = H$ as above. But \eqref{eq:lem:abelian variety} shows $\ker (s^* \circ \lambda)^\circ = \ker (i^* \circ \lambda)^\circ = H_\cdot^\perp$, where the latter equality is just the definition of $H_\cdot^\perp$. This shows $H_\cdot$ is the complement of $\ker (s^* \circ \lambda)^\circ$ as claimed.
		\end{proof}
		
		\begin{proof} [Proof of \Cref{prop:image of R_i}]
			To ease notation, let us drop the subscript $i$ throughout the proof. Hence $R$ means $R_i$, and the diagram
			\[\begin{tikzcd}[column sep=normal]
				R \arrow[r, "\forget"] \arrow[rd, "\delta"'] & A \times A \arrow[d, "-"] \\
				& A_\cdot
			\end{tikzcd}\]
			is the restriction of the original diagram to the irreducible component $R_i$.
			
			\medskip
			
			(1) The assertion is trivial when $\delta$ is surjective, so assume $\delta(R)$ has dimension $r < n$. Over a general point $a \in \delta(R) \subset A_\cdot$, the fiber $R_a = \delta^{-1}(a)$ has pure dimension $n-r$. Take a connected component $F \subset R_a$ and consider the restricted morphisms
			\[ s_F : F \to A ,\qquad t_F : F \to A \qquad\text{such that } \ t_F = t_a \circ s_F .\]
			Now \eqref{eq:sP = tP} induces the equality in $\Pic_F$:
			\[ (s_F^* \circ \lambda)(a) = s_F^* (t_a^* P \otimes P^{-1}) = t_F^* P \otimes s_F^* P^{-1} = \mathcal O_F .\]
			That is, we get $a \in \ker (s_F^* \circ \lambda)$ for every connected component $F \subset R_a$. Since $s_F$ is finite, \Cref{lem:abelian variety}(2) shows $\dim \ker (s_F^* \circ \lambda) \le n - (n-r) = r$. Varying the connected component $F \subset R_a$, the general point $a \in A_\cdot$ is contained in the subgroup
			\[ D_a \coloneq \bigcap_{F \subset R_a} \ker \Big( s_F^* \circ \lambda \Big) \qquad\text{of dimension } \le r .\]
			
			Therefore, there exists a Zariski open subset $U \subset \delta(R)$ that is contained in the union of proper subgroups
			\[ U \subset \bigcup_{a \in U} D_a ,\qquad D_a \subsetneq A_\cdot \text{ is a subgroup of dimension at most } r < n .\]
			But there exist only countably many proper subgroups of $A_\cdot$. This forces $U \subset D_a$ for a single $a \in U$, and therefore $\delta(R) \subset D_a$.\footnote{This argument is reminiscent of Kawamata's proof of his \cite[Theorem~13]{kaw81}.} Since any subgroup of $A_\cdot$ is of the form $D \times \pi_0$ for an abelian subvariety $D \subset A_\cdot$ and finite subgroup $\pi_0 \subset A_\cdot$, we have $\delta(R) \subset D + \tau$ for some torsion point $\tau \in A_\cdot$ and an abelian variety $D$. The sequence of dimension estimates
			\[ r = \dim \delta(R) \le \dim D \le r \]
			is now all equalities. This concludes $\delta(R) = D + \tau$ and $\dim D = r$.
			
			\medskip
			
			(2) By the equality assertion of \Cref{lem:abelian variety}(2), for every general point $a \in \delta(R) = D + \tau$, the fiber
			\[ \bar R_a = \bar \delta^{-1}(a) \cong s(R_a) \subset A \]
			is a finite union of abelian subtorsors of the form $K + \sigma$ for $K = D^\perp$. This shows that $\bar \delta^{-1} (U)$ is closed under the given $K$-action on $A \times A$. Since $\bar R$ is irreducible, such a $K$-action on its Zariski open subset uniquely extends to a $K$-action on the whole $\bar R$. Note that this automatically means $K$ acts fiberwise on $\bar \delta : \bar R \to A_\cdot$.
		\end{proof}

	\subsection{Composition of irreducible components}
		Assume further that $B^{sch}$ is a Fano variety, and let us begin the proof of \Cref{thm:delta(R_i)=A} following \S 11.4 of the paper.
		
		\medskip
		
		Use \Cref{prop:inertial component} to pick an inertial component $R_i \subset R$. If $\delta(R_i) = A_\cdot$ then we are done, so assume that (\Cref{cor:image of R_i})
		\begin{equation} \label{eq:D}
			D = \delta(R_i) \subset A_\cdot ,\qquad K = D^\perp \subset A_\cdot
		\end{equation}
		are abelian subvarieties of positive dimensions $r$ and $n-r$. The strategy will be to find a new inertial component $R_k \subset R$ whose difference $\delta(R_k)$ is strictly larger than $D$. To do so, in \Cref{lem:prop:inertial component 1}, we will first find a suitable non-identity component $R_j \subset R$ adjacent to the identity component ($R_0 \cap R_j \neq \emptyset$). In \Cref{lem:prop:inertial component 3}, the new inertial component $R_k$ will be obtained as the \emph{composition} of $R_i$ and $R_j$ using the groupoid structure of $R$.
		
		\medskip
		
		Let us first make the following notation explicit.
		
		\begin{definition}
			Let $K \subset A_\cdot$ be an abelian subvariety. For each closed subvariety $X \subset A$, we define $K+X \subset A$ as the image of $K \times X \subset A_\cdot \times A$ under the translation morphism $+ : A_\cdot \times A \to A$. Similarly, for $Z \subset A \times A$ we define $K+Z \subset A \times A$ under the diagonal $A_\cdot$-action on $A \times A$.
		\end{definition}
		
		\begin{lemma} [AI::Lemma~11.7] \label{lem:prop:inertial component 1}
			There exist a non-identity component $R_j$ together with a prime Weil divisor $E_j \subset R_j$ satisfying the following three conditions:
			\begin{enumerate}[label=\textnormal{(\alph*)}]
				\item $\delta(R_j) \subset A_\cdot$ is an abelian subvariety, and is not contained in $D$.
				\item $R_j$ intersects with $R_0$, and $E_j$ is an irreducible component of $R_0 \cap R_j$.\newline
				In particular, $s : R \to A$ induces an isomorphism $E_j \cong s(E_j) \subset A$.
				\item $K + s(E_j) = A$.
			\end{enumerate}
		\end{lemma}
		\begin{proof}
			Let us first construct a divisor $E_j$ in the identity component $R_0$. Let $T \subset A$ be the ramification divisor of $p^{sch} : A \to B \to B^{sch}$ arising in the Riemann--Hurwitz formula:
			\[ K_A (= \mathcal O_A) \sim_{\QQ} (p^{sch})^* K_{B^{sch}} \otimes \mathcal O_A(T) .\]
			Then $T$ is ample because $K_{B^{sch}}$ is anti-ample and $p^{sch}$ is finite. This in particular implies $K + T = A$ (otherwise we have $K + T = T$, violating the ampleness of $T$). Take $E_j \subset T$ to be any prime component such that (c) $K + E_j = A$. We may consider $E_j$ as a divisor in $R_0$ via an isomorphism $s : R_0 \to A$.
			
			\medskip
			
			Let us now find the component $R_j$. By \Cref{lem:codimension 1 description of R}, the morphism $s : R \to A$ is ramified along $E_j \subset R_0 \subset R$ with locally $e > 1$ number of branches. Take any non-identity component $R_j$ containing some of these branches. Then (b) $E_j \subset R_j \cap R_0$ is satisfied.
			
			\medskip
			
			It only remains to show (a) $\delta(R_j) \not \subset D$. \Cref{cor:image of R_i} shows $\delta(R_j) = D_j \subset A_\cdot$ is an abelian subvariety. Assume on the contrary that we have $D_j \subset D$. Then this implies the opposite inclusion $K_j = D_j^\perp \supset K$ and hence shows that $\bar R_j$ is closed under the $K$-action by \Cref{prop:image of R_i}(2). But $\bar R_0 = \Delta_A \subset A \times A$ is closed under the $K$-action as well, so their intersection $\bar R_j \cap \bar R_0$ is closed under the $K$-action. This means our irreducible component $E_j \subset R_0 = \bar R_0$ is also closed under the $K$-action, contradicting the condition (c) $K + E_j = A$.
		\end{proof}
		
		\begin{lemma} [AI::Lemma~11.8] \label{lem:prop:inertial component 2}
			Given $E_j \subset R_j$ in \Cref{lem:prop:inertial component 1}, there exists a prime Weil divisor $E_i \subset R_i$ satisfying the following three conditions:
			\begin{enumerate}[label=\textnormal{(\alph*)}]
				\item $\delta(E_i) = D$.
				\item $t(E_i) = s(E_j)$ for $t : R_i \to A$ and $s : R_j \to A$.
				\item $E_i$ contains a point of the form $(x,x,g)$ with $g \neq \id$.
			\end{enumerate}
		\end{lemma}
		\begin{proof}
			Since $R_i$ is inertial, we can take any point $e_0 = (x_0,x_0,g_0) \in R_i$ with $g_0 \neq \id$. Our first goal is to find another point
			\[ e_1 = (x_1, x_1, g_1) \in R_i \qquad\text{with}\quad x_1 \in s(E_j) \ \text{ and }\ g_1 \neq \id .\]
			By the condition $K + s(E_j) = A$ in \Cref{lem:prop:inertial component 1}(c), there exists $k \in K$ such that $x_1 \coloneq k + x_0 \in s(E_j)$. \Cref{prop:image of R_i}(2) constructed a $K$-action on $\bar R_i$, and this action lifts to a $K$-action on its normalization $Z_i$ by \Cref{lem:bar R_i}. Take any lift $\tilde e_0 \in Z_i$ and its $K$-orbit $K \cdot \tilde e_0 \subset Z_i$. Its image in $R_i$ passes through both $e_0 = (x_0, x_0, g_0)$ and another point of the form $e_1 = (x_1, x_1, g_1)$. Now $g_0 \neq \id$ implies $g_1 \neq \id$ by \Cref{lem:I_A}.
			
			\medskip
			
			Since $t : R_i \to A$ is finite surjective (\Cref{lem:R}), we may take any irreducible component $E_i \subset t^{-1}(s(E_j))$ passing through this point $e_1 = (x_1, x_1, g_1)$. It is a prime divisor in $R_i$ with $t(E_i) = s(E_j)$. This achieves both (b) and (c).
			
			\medskip
			
			It remains to show (a) $\delta(E_i) = D$. Recall the diagram
			\[\begin{tikzcd}[column sep=normal]
				R_i \arrow[r, "\forget"] \arrow[rd, "\delta"'] & \bar R_i \arrow[d, "\bar \delta"] \\
				& A_\cdot
			\end{tikzcd}\]
			and write $\bar E_i = \forget(E_i) \subset \bar R_i$. Let us first show $K + \bar E_i = \bar R_i$ under the $K$-action on $\bar R_i$ in \Cref{prop:image of R_i}(2). Taking projection of its left hand side is
			\[ p_2 (K + \bar E_i) = K + t(E_i) = K + s(E_j) = A \qquad\text{by (b) and \Cref{lem:prop:inertial component 1}(c)} .\]
			Hence $K + \bar E_i \subset \bar R_i$ dominates $A$ via the finite surjective morphism $p_2 : \bar R_i \to A$. But $\bar R_i$ is irreducible, so this is enough to show $K + \bar E_i = \bar R_i$. This claim implies the sequence of equalities
			\[ \delta(E_i) = \bar \delta(\bar E_i) = \bar \delta(K + \bar E_i) = \bar \delta (\bar R_i) = \delta(R_i) = D .\]
			Here the second equality followed because $\bar \delta$ is invariant under the $K$-action on $\bar R_i$.
		\end{proof}
		
		\begin{lemma} \label{lem:prop:inertial component 3}
			Under \Cref{lem:prop:inertial component 1}--\ref{lem:prop:inertial component 2}, there exists a new inertial component $R_k \subset R$ satisfying $\delta(R_k) \supsetneq D$.
		\end{lemma}
		\begin{proof}
			Using a morphism $t : E_i \to s(E_j)$ in \Cref{lem:prop:inertial component 2}(b) and an isomorphism $s : E_j \to s(E_j)$ in \Cref{lem:prop:inertial component 1}(b), we can construct the graph morphism $\iota : E_i \hookrightarrow R_{ij}$ by $\iota(e) = (e, s^{-1}(t(e)))$ as in the following diagram:
			\[\begin{tikzcd}
				& E_j \arrow[d, "s"', "\cong"] \\
				E_i \arrow[r, "t"] & s(E_j)
			\end{tikzcd},\qquad\quad\begin{tikzcd}[column sep=normal]
				E_i \arrow[r, hook, "\iota"] \arrow[rd, "\id"'] \arrow[rr, bend left=35, "s^{-1} \circ t"] & R_{ij} \arrow[r, "p_2"] \arrow[d, "p_1"] & R_j \arrow[d, "s"] \\
				& R_i \arrow[r, "t"] & A
			\end{tikzcd}.\]
			Take any irreducible component $T \subset R_{ij}$ that contains $\iota(E_i)$. Then the image
			\[ R_k = c(T) \subset R \]
			under the composition map $c : R \times_A R \to R$ is an irreducible component by \Cref{lem:irreducible components of R^2}.
			
			\medskip
			
			Let us first show that $R_k$ is inertial. Recalling the point $e_0 = (x_0,x_0,g_0) \in E_i$ with $g_0 \neq \id$ in \Cref{lem:prop:inertial component 2}(c), we have
			\[ \iota(e_0) = (e_0, s^{-1}(t(e_0))) = ((x_0,x_0,g_0), (x_0,x_0,\id)) \in T .\]
			Their composition is $c((x_0,x_0,g_0), (x_0,x_0,\id)) = (x_0,x_0,g_0) = e_0 \in R_k$, so $R_k$ an inertial component (and shows $R_i \cap R_k \neq \emptyset$).
			
			\medskip
			
			It remains to show $\delta(R_k) \supsetneq D$. Notice first $\delta(R_k) = (\delta \circ c)(T)$. Applying the relation $\delta \circ c = \delta \circ p_1 + \delta \circ p_2$ in \Cref{prop:composition} to $\iota(e) \in T$ for $e \in E_i$ yields
			\[ (\delta \circ c)(\iota(e)) = (\delta \circ p_1)(\iota(e)) + (\delta \circ p_2)(\iota(e)) = \delta(e) + \delta(s^{-1}(t(e))) = \delta(e) .\]
			Here the last equality follows from $\delta \equiv 0$ on $E_j \subset R_j$ by \Cref{lem:prop:inertial component 1}(b). Since $\delta(E_i) = D$ by \Cref{lem:prop:inertial component 2}(a), we get an inclusion $(\delta \circ c)(T) \supset D$.
			
			On the other hand, pick a point $t \in T$ with $(\delta \circ p_2)(t) \notin D$ using the surjectivity of $p_2 : T \to R_j$ by \Cref{lem:irreducible components of R^2} and $\delta(R_j) \not \subset D$ by \Cref{lem:prop:inertial component 1}(a). Since $(\delta \circ p_1)(t) \in \delta(R_i) = D$, this shows
			\[\begin{aligned}
				(\delta \circ c)(t) = (\delta \circ p_1)(t) + (\delta \circ p_2)(t) \ \ \notin \ D
			\end{aligned},\]
			so the inclusion is strict as desired.
		\end{proof}
		
		\begin{proof} [Proof of \Cref{thm:delta(R_i)=A}]
			If $B$ is not a scheme, then we can find an inertial component $R_i \subset R$. If $\delta(R_i) = D$ is the entire $A_\cdot$, then we are done. If not, then use \Cref{lem:prop:inertial component 1}--\ref{lem:prop:inertial component 3} to find another inertial component $R_k \subset R$ with $\delta(R_k) \supsetneq D$. Induction completes the proof.
		\end{proof}

\section{Derivative of the presentation}
	This section summarizes AI::Lemma~11.4--11.5.
	
	\begin{theorem} [AI::Lemma~11.5] \label{thm:delta(R_i) is not A}
		Let $p : A \to B$ be as usual and assume \Cref{thm:main}(ii). Then every inertial component $R_i \subset R$ satisfies
		\[ \delta(R_i) \subsetneq A .\]
	\end{theorem}
	
	The theorem does not assume $B^{sch}$ is projective; we do not need \Cref{prop:image of R_i}. In fact, the proof of this theorem is almost purely local on $A$. Besides developing a general theory, the two main ingredients are (1) interpretation of the top Chern class of a vector bundle as the vanishing locus of its regular section, and (2) translating the information on the difference morphism $\delta$ to a condition on the derivative $dp : T_A \to p^*T_B$ of the morphism $p$.

	\subsection{Incidence correspondence and the Euler characteristic assumption}
		At each point $x \in A$, the morphism $p : A \to B$ induces a derivative homomorphism
		\[ dp_x : T_x A \to T_{p(x)} B .\]
		Again for each automorphism $g$ of $p(x)$, the tangent space $T_{p(x)} B$ admits the $g$-eigenspace decomposition $T_{p(x)}^g B \oplus (\text{complement})$. Composing $dp_x$ with the projection $T_{p(x)} B \to T_{p(x)}^g B$ yields a linear map
		\begin{equation} \label{eq:alpha}
		\begin{tikzcd}[sep=normal]
			\alpha_{(x,g)} : T_x A \arrow[r, "dp_x"] & T_{p(x)} B \arrow[r, twoheadrightarrow] & T_{p(x)}^g B
		\end{tikzcd}.
		\end{equation}
		In this section, we prove the following result about this map \eqref{eq:alpha}.
		
		\begin{proposition}[AI::Lemma~11.4] \label{prop:vector 1}
			Assume \Cref{thm:main}(ii). Then a general vector field $v \in H^0 (A, T_A)$ satisfies $\alpha_{(x,g)} (v) \neq 0$ for every $x \in A$ and every automorphism $g$ of $p(x) \in B$.
%			\[ \alpha_{(x,g)} (v_x) \neq 0 \quad\text{for every $x \in A$ and every automorphism $g$ of $p(x) \in B$} .\]
%			satisfies $(\id + g + \cdots + g^{d-1})(\xi) \neq 0$ for every automorphism $g$ of $p(x)$ of order $d$.
		\end{proposition}
		
		To understand this statement more systematically, recall the cartesian diagram between the inertia stack and stabilizer scheme
		\[\begin{tikzcd}
			I_A \arrow[r, "p"] \arrow[d, "j"] & I_B \arrow[d, "j"] \\
			A \arrow[r, "p"] & B
		\end{tikzcd}.\]
		The pullback $p^* \Omega_B \to \Omega_A$ models the dual of the derivative $dp_x$ above. Its restriction to the splitting injection $\Omega_{I_B} \hookrightarrow j^* \Omega_B$ is a homomorphism of vector bundles on $I_A$ modeling the dual of \eqref{eq:alpha}\footnote{This is the dual of the homomorphism $\alpha : j^* T_A \to p^* T_{I_B}$ defined after AI::(11.6) in p.89.}
		\[\begin{tikzcd}[sep=normal]
			\check \alpha : p^* \Omega_{I_B} \arrow[r, hook] & p^* j^* \Omega_B = j^* p^* \Omega_B \arrow[r] & j^* \Omega_A \ (\cong \mathcal O_{I_A}^{\oplus n})
		\end{tikzcd}.\]
		\Cref{prop:vector 1} studies the property of the set
		\[ L = \{ (x,g,v) \in I_A \times \AA^n : \alpha_{(x,g)}(v) = 0 \} \]
		parametrizing the \emph{fiberwise} kernel of $\alpha$. The set admits a scheme structure as follows.
		
		\begin{lemma} \label{lem:prop:vector 1 step 1}
			$L$ is a scheme with an incidence correspondence
			\[\begin{tikzcd}
				L \arrow[r, "q"] \arrow[d] & \AA^n \\
				I_A
			\end{tikzcd}.\]
		\end{lemma}
		\begin{proof}
			Take the cokernel of $\check \alpha : p^* \Omega_{I_B} \to \mathcal O_{I_A}^{\oplus n}$ above. It admits a surjection $\mathcal O_{I_A}^{\oplus n} \twoheadrightarrow \coker \check \alpha$. We can define $L$ as the linear scheme associated to $\coker \check \alpha$:
			\[ L \coloneq \Spec_{I_A} \Big( \Sym^* (\coker \check \alpha) \Big) \ \ \subset \ \ I_A \times \AA^n = \Spec_{I_A} \Big( \Sym^* (\mathcal O_{I_A}^{\oplus n}) \Big) .\]
			This $L$ is the \emph{fppf dual} of $\coker \check \alpha$, so it parametrizes the fiberwise kernel $\ker \alpha_{(x,g)}$ for every point $(x,g) \in I_A$.
		\end{proof}
		
		\begin{lemma} \label{lem:prop:vector 1 step 2}
			The degeneracy locus
			\[ I_{A, \le r} \coloneq \{ (x,g) \in I_A : \rk \alpha_{(x,g)} \le r \} \]
			is a closed subscheme of $I_A$ of dimension $\le r$.
		\end{lemma}
		\begin{proof}
			Note that $\alpha_{(x,g)}$ and $\check \alpha_{(x,g)}$ have the same rank, so this set is the degeneracy locus of a vector bundle homomorphism $\check \alpha$, see \cite[\S 14]{ful:intersection}. Hence it is a closed subscheme of $I_A$.
			
			To show the dimension claim, we may show instead: $\check \alpha$ has a generic rank $\ge s$ along every irreducible subvariety $i : Y \hookrightarrow I_A$ of dimension $s$.
			\[\begin{tikzcd}[column sep=normal]
				Y \arrow[r, hook, "i"] \arrow[rd] & I_A \arrow[r, "p"] \arrow[d, "j"] & I_B \\
				& A
			\end{tikzcd}\]
			Notice that $j \circ i$ is unramified because $j$ is unramified (\Cref{lem:I_A}). Hence we have a surjective homomorphism $(j \circ i)^* \Omega_A \twoheadrightarrow \Omega_Y$. The composition
			\[\begin{tikzcd}[sep=normal]
				(p \circ i)^* \Omega_{I_B} \arrow[r, "\check \alpha_{|Y}"] & (j \circ i)^* \Omega_A \arrow[r, twoheadrightarrow] & \Omega_Y
			\end{tikzcd}\]
			is nothing but the cotangent sequence of the morphism $p \circ i$. This homomorphism is generically surjective because the finite morphism $p \circ i$ is generically unramified. Therefore, it has a generic rank $s$. Hence $\check \alpha_{|Y}$ has a generic rank $\ge s$ by the surjectivity of the second homomorphism.
		\end{proof}
		
		\begin{corollary} \label{lem:prop:vector 1 step 3}
			\begin{enumerate}
				\item The set $I_{A, 0} = \{ (x,g) \in I_A : \alpha_{(x,g)} \equiv 0 \}$ is finite.
				\item Every general $v \in H^0 (A, T_A)$ satisfies $\alpha_{(x,g)} (v) \neq 0$ for every $(x,g) \in I_{A, \ge 1}$.
			\end{enumerate}
		\end{corollary}
		\begin{proof}
			(1) Immediate from \Cref{lem:prop:vector 1 step 2}. (2) The morphism $L \to I_A$ is an $\AA^{n-r}$-bundle over $I_{A, r}$, which has dimension $\le r$ by \Cref{lem:prop:vector 1 step 2}. Hence every irreducible component of $L$ has dimension $\le n$.
			
			Let $L_i \subset L$ be an irreducible component, and consider the projection morphism in \Cref{lem:prop:vector 1 step 1}:
			\[ q : L \to \AA^n ,\qquad (x,g,v) \mapsto v .\]
			If $q(L_i) = \AA^n$, then this forces $L_i$ to be of the form
			\[ L_i = \{ (x,g) \} \times \AA^n \qquad\text{for}\quad (x,g) \in I_{A, 0} .\]
			If $q(L_i) \subsetneq \AA^n$, then $q(L_i)$ a proper constructible subset of $\AA^n$ by Chevalley's theorem. The union of all such proper subsets $q(L_i)$ is still proper, so its complement contains a Zariski open subset $U \subset \AA^n$. Combining these two, we obtain the description
			\[ q^{-1}(v) = I_{A, 0} \times \{ v \} \qquad\text{for every } v \in U \subset \AA^n ,\]
			which is the desired claim.
		\end{proof}
		
		\begin{proof}[Proof of \Cref{prop:vector 1}]
			All the above lemmas and corollary were completely formal and the assumption in \Cref{thm:main}(ii) was never used. Now we prove that the given assumption implies $I_{A,0} = \emptyset$ and completes the proof of the proposition.
			
			\medskip
			
			Take the Zariski open subset $U \subset \AA^n$ in \Cref{lem:prop:vector 1 step 3} and let $v \in U$ be any generic vector. We have shown $q^{-1}(v) = I_{A, 0} \times \{ v \}$, which means that the global section $\alpha(v)$ of $p^* T_{I_B}$ has a finite vanishing set $I_{A, 0}$. Therefore, according to the geometric interpretation of the top Chern class of vector bundles \cite[Proposition~14.1]{ful:intersection} \cite[Tag~0FA9]{stacks-project}, we can show $I_{A, 0} = \emptyset$ once the top Chern class of $p^* T_{I_B}$ vanishes. The assumption on the Euler characteristic of each connected component means the top Chern class of $T_{I_B}$ vanishes already on $I_B$ by the Gauss--Bonnet theorem in \cite[Proposition~1.6]{beh09}. This completes the proof.
		\end{proof}

	\subsection{Local structure of the equivalence relation}
		\Cref{thm:delta(R_i) is not A} is essentially a \emph{local} study of the morphism $p : A \to B$. The idea is to translate the surjectivity of the difference map $\delta$ into a local property of $p$. This is the content of the following observation, which is elementary but crucial.
		
		\begin{lemma} \label{lem:identification of tangent spaces}
			Let $\Delta$ be a smooth affine curve with a reference point $0$, and $c : \Delta \to A \times A$ a nonconstant morphism with $c(0) = (x_0, y_0)$. Consider three compositions
			\[ c_1 \coloneq p_1 \circ c ,\quad c_- \coloneq - \circ c ,\quad c_2 \coloneq p_2 \circ c \qquad\text{for}\quad \begin{tikzcd}
				\Delta \arrow[r, "c"] & A \times A \arrow[d, "p_1"'] \arrow[rd, "-"'] \arrow[r, "p_2"] & A \\
				& A & A_\cdot
			\end{tikzcd}.\]
			Then we have an equality
			\[ dc_- (\partial_t) = dc_2(\partial_t) - dc_1(\partial_t) \qquad\text{for}\quad \partial_t = \tfrac{\partial}{\partial t} \in T_0 \Delta \]
			under the identifications $T_{x_0} A = T_{y_0 - x_0} A_\cdot = T_{y_0} A = H^0 (A, T_A)$.
		\end{lemma}
		\begin{proof}
			The easiest way to see this is by setting $A = \CC^n/\Gamma$. The statement is analytically local, so we can pretend $\Delta$ is a unit disc, $A = T_0 A = \CC^n$, and $x_0 = y_0 = 0$. Then the subtraction map on $A$ coincides with that on $\CC^n$. Now the lemma is saying nothing but $(f-g)'(0) = f'(0) - g'(0)$ for two holomorphic functions $f, g : \CC \to \CC^n$.
		\end{proof}
		
		To prove \Cref{thm:delta(R_i) is not A}, assume on the contrary that $R_i \subset R$ is an inertial component with surjective $\delta : R_i \to A_\cdot$. Then \Cref{lem:identification of tangent spaces} has the following consequence.
		
		\begin{proposition} [AI::Lemma~11.5] \label{prop:vector 2}
			Let $R_i$ be an inertial component with surjective $\delta : R_i \to A_\cdot$. Then for every general $v \in H^0 (A, T_A)$, there exist a point $(x,x,g) \in R_i$ with $g \neq \id$ and a local presentation $f : (U, x) \to (V, 0)$ of $p$ such that
			\[ df_x(v) \in \im (\id - g) \quad (\subset T_0 V) .\]
		\end{proposition}
		\begin{proof}
			Every tangent vector $v \in H^0 (A, T_A) = T_0 A_\cdot$ is realized as the direction vector $dc_- (\partial_t)$ of a morphism $c_- : \Delta \to A_\cdot$ from a smooth affine curve $\Delta$ with $c_-(0) = 0$. The surjectivity assumption on $\delta : R_i \to A_\cdot$ means that $\delta$ is generically \'etale. Thus, a general choice of $v \in T_0 A_\cdot$ makes $\delta$ generically \'etale over the curve $\Delta$. We will shrink $\Delta$ around $0$ whenever needed.
			
			Pick a point $(x_0, x_0, g_0) \in R_i$ with $g_0 \neq \id$, and take a connected component $F \subset R_i \cap I_A$ containing this point. Since $\delta : R_i \to A_\cdot$ is generically \'etale over $\Delta$, we can lift this curve to a rational section $c : \Delta \dashrightarrow R_i$ through $F$. The valuative criterion ($\delta$ is proper) fills in the missing value $c(0)$ and constructs a lifting
			\[\begin{tikzcd}[column sep=normal]
				& R_i \arrow[d, "\delta"] \\
				\Delta \arrow[ru, bend left=25, "c"] \arrow[r, "c_-"] & A_\cdot
			\end{tikzcd}.\]
			By construction, $c(0)$ passes through the connected component $F$ of $R_i \cap I_A$, which contains the given point $(x_0, x_0, g_0)$ with $g_0 \neq \id$. Recalling \Cref{lem:I_A}, the point $c(0)$ must be of the form $(x_1, x_1, g_1)$ with $g_1 \neq \id$.
			
			\medskip
			
			Apply \Cref{lem:neighborhood of inertial component} to the point $(x_1, x_1, g_1)$ to take a local presentation $f : U \to V$ of $p : A \to B$ and a local neighborhood $W$ of $R_i$. Taking its image $\bar W = \forget(W) \subset U \times U$, we have the equality
			\begin{equation} \label{eq:f(y) = g f(x)}
				f(y) = g_1 \cdot f(x) \qquad\text{for}\quad (x,y) \in \bar W
			\end{equation}
			Denote the curve on $\bar W$ by $c : \Delta \to \bar W$ for simplicity, and write $c = (c_1, c_2)$ for $c_1, c_2 : \Delta \to U$.
			\[\begin{tikzcd}[column sep=normal]
				\Delta \arrow[r, "c"] \arrow[rd, bend right=25, "c_1"'] \arrow[rr, bend left=35, "c_2"] & \bar W \arrow[r, "p_2"] \arrow[d, "p_1"] & U \\
				& U
			\end{tikzcd}\]
			\Cref{lem:identification of tangent spaces} claims $v = dc_2(\partial_t) - dc_1(\partial_t)$, so applying $df$ and using \eqref{eq:f(y) = g f(x)} shows the desired claim
			\begin{align*}
				df(v) &= d(f \circ c_2)(\partial_t) - d(f \circ c_1)(\partial_t) \\
				&= g \cdot d(f \circ c_1)(\partial_t) -  d(f \circ c_1)(\partial_t) = (g - \id) \cdot d(f \circ c_1)(\partial_t) .\qedhere
			\end{align*}
		\end{proof}
		
		\begin{proof}[Proof of \Cref{thm:delta(R_i) is not A}]
			The two conclusions \Cref{prop:vector 1} and \Cref{prop:vector 2} are incompatible. To see this, recall that $\alpha_{(x,g)}(v)$ was defined to be $dp(v) = df_x(v)$ projected down to the $g$-invariant subspace of $T_{p(x)} V = T_0 V$. The projection homomorphism is explicitly described as the averaging map
			\[ T_0 V \twoheadrightarrow T_0^g V ,\qquad v \mapsto \tfrac{1}{d}(\id + g + \cdots + g^{d-1})(v) ,\]
			where $d$ is the order of the linear map $g$ on $T_0V$. Hence \Cref{prop:vector 1} says that a general $v \in H^0 (A, T_A)$ satisfies $(\id + g + \cdots + g^{d-1})(df_x(v)) \neq 0$ for every $(x,g) \in I_A$.
			
			However, \Cref{prop:vector 2} constructed $(x_0, g_0) \in I_A$ such that $df_{x_0}(v) = (\id - g_0)(\xi)$ for some $\xi \in T_0 V$. This computes
			\[ (\id + g_0 + \cdots + g_0^{d-1})(df_{x_0}(v)) = (\id + g_0 + \cdots + g_0^{d-1})(\id - g_0)(\xi) = (\id - g_0^d)(\xi) = 0 ,\]
			which is the desired contradiction.
		\end{proof}

\section{Proof of main theorems} \label{sec:proof of main theorems}
	The proof of \Cref{thm:main} is now immediate from the two main results of the previous sections.
	
	\begin{proof} [Proof of \Cref{thm:main}]
		If $B$ is not a scheme, then (i) implies \Cref{thm:delta(R_i)=A} and (ii) implies \Cref{thm:delta(R_i) is not A}. These are incompatible.
	\end{proof}
	
	The following proposition reduces \Cref{thm:base} to \Cref{thm:main}.
	
	\begin{proposition} \label{prop:reduction of theorem}
		Let $X$ be a projective hyper-K\"ahler manifold and $f : X \to B$ a Lagrangian fibration to a smooth proper DM stack $B$ with trivial stabilizers in codimension $1$. Then
		\begin{enumerate}
			\item $B^{sch}$ is a Fano variety.
			\item Every (non-identity) connected component of $I_B$ has Euler characteristic $0$.
		\end{enumerate}
	\end{proposition}
	
	Again, the \emph{Euler characteristic} of a DM stack here means its stacky-Euler characteristic, which is a rational number in general as in \cite[\S 1.3]{beh09}. The first item is a well-known result of Matsushita. The second item is new; OpenAI's idea was to transform another well-known theorem of Matsushita into a K-theoretic form. The rest of this note proves \Cref{prop:reduction of theorem} and concludes.

	\subsection{Euler characteristic of the inertia stack}
		Let us start with some formal discussions on vector bundles on smooth, separated, and finite type DM stacks $S$. Given a vector bundle $F$ of rank $r$ on $S$, we define\footnote{Our notation is not standard. OpenAI followed the standard notation $\lambda_{-1}(F)$ in the context of $\lambda$-rings.}
		\[ \wedge^* (F) = \sum_{i=0}^r (-1)^i \cdot \left[ {\textstyle \bigwedge^i} F \right] \qquad\text{in}\quad K^0 (S) .\]
		Similarly, one can define $\wedge^* (\alpha)$ for any class $\alpha \in K^0 (S)$.
		
		\begin{lemma} \label{lem:reduction of theorem:Chern class}
			Let $S$ be a smooth separated DM stack and $\ch : K^0(S)_{\bar \QQ} \to \CH^* (S)_{\bar \QQ}$ the Chern character homomorphism. Then
			\begin{enumerate}
				\item $\wedge^* (\alpha + \beta) = \wedge^*(\alpha) \cdot \wedge^* (\beta)$ for $\alpha, \beta \in K^0 (S)$.
				\item For each rank $r$ vector bundle $F$ on $S$, we have
				\[\begin{aligned}
					\ch \Big( \wedge^*(F) \Big) &= (-1)^r \cdot c_r (F) \ + (\text{terms of degree $> r$}) ,\\
					\ch_0 \Big( \wedge^* (c \cdot [F]) \Big) &= (1 - c)^r \qquad\text{for}\quad c \in \bar \QQ .
				\end{aligned}\]
			\end{enumerate}
		\end{lemma}
		\begin{proof}
			(1) This is standard. (2) When $F$ is a line bundle, both are easily deduced. When $F$ is a vector bundle of rank $r > 1$, a standard method of Grothendieck (passing to the flag variety-bundle over $S$) reduces to the case when $F$ is a direct sum of line bundles. Now use (1) and the fact that $\ch$ is a ring homomorphism.
		\end{proof}
		
		Let us return back to \Cref{prop:reduction of theorem} and assume $B$ is the smooth DM stack base of a Lagrangian fibration. Our references for some facts about inertia stacks are \cite[Theorem~3.15]{toen99} and \cite[\S \S 2.1--2.2]{tseng10}.
		
		\begin{lemma} [AI::Lemma~11.3] \label{lem:reduction of theorem:K0}
			Every rank $r$ vector bundle $E$ on $B$ satisfies an equality
			\[ ([E] - r)^m \cdot \wedge^* (\Omega_B) = 0 \qquad\text{for } \ m \gg 0 \ \text{ in } K^0(B) .\]
		\end{lemma}
		\begin{proof}
			Every proper flat Lagrangian fibration $f : X \to B$ from a quasi-projective symplectic variety $X$ satisfies the equality $R^if_* \mathcal O_X = \Omega_B^i$ by \cite{mat05, sch23}. This means $[Rf_* \mathcal O_X] = \wedge^*(\Omega_B)$, or similarly using the projection formula
			\[ Rf_*([f^*E]^i) = [E]^i \cdot \wedge^*(\Omega_B) \qquad\text{in}\quad K^0(B) .\]
			Now the rank $r$ vector bundle $f^* E$ on a smooth variety $X$ satisfies the identity (because in this case, the Chern character homomorphism is an isomorphism)
			\[ ([f^* E] - r)^m = 0 \qquad\text{for } \ m > 2n \quad \text{in } K^0(X) .\]
			Applying $Rf_*$ to this identity gives the desired claim.
		\end{proof}
		
		Let $S \subset I_B$ be a non-identity connected component of the inertia stack, and write $j : S \subset I_B \to B$ for the composition. Let $E$ be a vector bundle on $B$ and $j^* E$ its pullback to $S$. Locally around a point $(z,g) \in S$, the vector bundle $j^*E$ is $g$-equivariant and hence admits the eigenbundle decomposition $j^* E = \bigoplus_{i=0}^{k-1} E^{(i)}$ where $k$ is the order of $g$, $E^{(i)}$ is the subbundle with eigenvalue $\zeta^i$, and $\zeta$ is a primitive $k$-th root of unity. Such a decomposition patches together into
		\[ j^* E = \bigoplus_{i=0}^{k-1} E^{(i)} \qquad\text{on}\quad S .\]
		Following \cite[Theorem~3.15]{toen99}, we define a ring homomorphism $\phi$ by
		\begin{equation} \label{eq:Toen}
			\phi : K^0(B) \to K^0(S)_{\bar \QQ} ,\qquad [E] \mapsto \sum_{i=0}^{k-1} \zeta^i \cdot [E^{(i)}] .
		\end{equation}
		Composing it with the Chern character map of $S$ is the \emph{inertial Chern character} homomorphism
		\[ \ch \circ \phi : K^0(B) \to \CH^*(S)_{\bar \QQ} .\]
		
		\begin{proof}[Proof of \Cref{prop:reduction of theorem}]
			(1) This is a well-known result of Matsushita \cite[Theorem~2]{mat99}, which uses the global hyper-K\"haler geometry of $X$. (2) Let $S$ be a non-identity component of $I_B$, say of dimension $d < n$, with a morphism $j : S \to B$. Say the stabilizers in $S$ have order $k$ and $\zeta$ is a primitive $k$-th root of unity as above.
			
			\medskip
			
			Start with the equality in \Cref{lem:reduction of theorem:K0} (for an arbitrary vector bundle $E$ on $B$). Applying the homomorphism $\phi$ \eqref{eq:Toen} to this equality and using \Cref{lem:reduction of theorem:Chern class}(1), we obtain
			\[ \Big( \phi(E) - r \Big)^m \cdot \prod_{i=0}^{k-1} \Big( \wedge^* (\zeta^i \cdot [\Omega_B^{(i)}]) \Big) = 0 \qquad\text{in}\quad K^0 (S)_{\bar \QQ} .\]
			Now subsequently apply the Chern character to it. \Cref{lem:reduction of theorem:Chern class}(2) computes
			\[\begin{aligned}
				\ch \Big( \wedge^* (\Omega_B^{(0)}) \Big) &= \ch \Big( \wedge^* (\Omega_S) \Big) = (-1)^d \cdot c_d (\Omega_S) = c_d (T_S) ,\\
				\ch_0 \Big( \wedge^* (\zeta^i \cdot [\Omega_B^{(i)}]) \Big) &\neq 0 \qquad\text{for}\quad 0 < i < k .
			\end{aligned}\]
			Furthermore, we have a direct computation
			\[ \ch_0 \Big( \phi(E) - r \Big) = \sum_{i=0}^{k-1} r_i \cdot \zeta^i - r \qquad\text{for}\quad r_i = \rk E^{(i)} .\]
			Because $\sum_i r_i = r$ is the rank of $E$, this Chern class is nonzero iff $r_0 < r$. Simply take a vector bundle $E$ on $B$ satisfying this property; for example, the paper used $E = T_B$ in AI::Lemma~11.3--11.4 because $T_S = T_B^{(0)}$ is a proper subbundle of $T_B$ (same works for $E = \Omega_B$).
			
			\medskip
			
			Combining all the above computations yields the equality AI::(11.6)
			\[ c_d(T_S) = 0 \qquad\text{in}\quad \CH^*(S)_{\bar \QQ} .\]
			Our argument so far never used the properness of $B$ and holds locally on $B$. When $B$ (and hence $S$) is proper, a prettier way of saying this is $\chi(S) = 0$. This completes the proof.
		\end{proof}
		
		\Cref{prop:reduction of theorem}(2) implies that \emph{every} connected component $I_B$ has dimension $\ge 1$.\footnote{See right after AI::(11.6) in p.89.} Consequently, $B^{sch}$ has no isolated singularities. In fact, such a result was already proved by M\"uller--Xu in \cite{mul-xu26}. Their paper precedes OpenAI's paper and the two proofs have similarities. We invite the interested readers to check M\"uller--Xu's paper, which may have motivated OpenAI's approach to \Cref{prop:reduction of theorem}(2).
		
		\begin{theorem} [{\cite[Theorem~1.3]{mul-xu26}}]
			Let $f : X \to B$ be a proper Lagrangian fibration of a smooth quasi-projective symplectic variety $X$ to a smooth DM stack $B$. Then $I_B$ has no components with dimension $0$. In particular, $B^{sch}$ has no isolated singularities.
		\end{theorem}
		\begin{proof}
			The main result of their paper claims that $B$ has no isolated singularities, but in fact the same proof concludes this slightly stronger claim. Let us quickly review their method to see this.
			
			\medskip
			
			The statement is local on $B$, so we may assume $B = [V/G]$ where $V$ is a smooth variety with a $G$-action fixing $0 \in V$. Taking the base change by $V \to B$ (and writing $B$ for $V$ for simplicity) yields a $G$-equivariant Lagrangian fibration $f : X \to B$ where the $G$-action on $X$ is free, and that on $B$ has a fixed point $0 \in B$. Denote by $F = X_0$ its central fiber, which admits a free $G$-action. Matsushita's isomorphism $R^if_* \mathcal O_X = \Omega_B^i$ is interpreted as an identification $H^i (F, \mathcal O_F) = \wedge^i (T_0^* B)$, which is also $G$-equivariant in this case.
			
			\medskip
			
			Let $g (\neq \id) \in G$ be any non-identity element. Then $g$ acts on $F$ freely. A version of the holomorphic Lefschetz fixed point formula in \cite{baum-ful-quart79} (see also \cite[Theorem~2.5]{mul-xu26}) says that the $g$-action on $H^* (F, \mathcal O_F) = \wedge^* (T_0^*B)$ satisfies $\det (g - \id) = 0$. This means that the $g$-action on $T_0^* B$ has a nontrivial $g$-eigenspace with eigenvalue $1$, or equivalently the $g$-action on $B$ has a fixed locus $B^g$ of dimension $\ge 1$. This means every connected component of $I_B$ has dimension $\ge 1$.
		\end{proof}
		
		M\"uller--Xu presented some examples of Lagrangian fibrations $f : X \to B$ (with non-proper $B$) where $I_B$ has components of dimension $\ge 1$. As shown in the proof of \Cref{prop:reduction of theorem}(2), their tangent bundles must have trivial top Chern classes.

\printbibliography
\end{document}